\documentstyle[12pt]{article}
\input amssym.def
 \input amssym.tex

\def\vv{{\underline{v}}}
\def\tt{{\underline{t}}}

\def\P{\Bbb P}

\def\L{\Bbb L}
\def\Z{\Bbb Z}

\def\C{\Bbb C}

\newtheorem{theorem}{Theorem}

\newenvironment{definition}
{\smallskip\noindent{\bf Definition\/}:}{\smallskip\par}

\newenvironment{remark}
{\smallskip\noindent{\bf Remark\/}.}{\smallskip\par}

\title{Integration with respect to Euler characteristic over
the projectivization of the space of functions and the Alexander
polynomial of a plane curve singularity.}
\author{A.Campillo
\and F.Delgado \thanks{First two authors were partially supported by
DGICYT PB97-0471 and by Junta de Castilla y Le\'on:
VA51/97. Address:
University of Valladolid, Dept. of Algebra, Geometry and Topology,
47005 Valladolid, Spain.
E-mail: campillo\symbol{'100}cpd.uva.es, fdelgado\symbol{'100}agt.uva.es}
\and S.M.Gusein--Zade \thanks{Partially supported by grants
RFBR--98--01--00612 and INTAS--97--1644. During the work on the paper
the author enjoyed the hospitality of the University of Nice.
Address: Moscow State University,
Dept. of Mathematics and Mechanics, Moscow, 119899, Russia.
E-mail: sabir\symbol{'100}mccme.ru}}
\date{}
\begin{document}

\def\Aut{{\operatorname{Aut}}}
\def\eps{\varepsilon}

\maketitle

For a reduced plane curve singularity $C=\bigcup\limits_{i=1}^{r}C_i$ 
($C_i$ are its irreducible components), let $\Delta^C(t_1, \ldots, t_r)$
be the Alexander polynomial of the link
$C\cap S_\varepsilon^3\subset S_\varepsilon^3$
for $\varepsilon > 0$ small enough (see, e.g., {\cite{EN}}).
We fix the Alexander polynomial $\Delta^C(t_1, \ldots, t_r)$
(in general it is defined only up to multiplication by a monomial
$\pm t_1^{n_1}\cdot\ldots\cdot t_r^{n_r}$) assuming that it is really
a polynomial (i.e., it does not contain variables in negative powers)
and that $\Delta^C(0, \ldots, 0) = 1$.
Let $\zeta_C(t)$ be the zeta-function of the classical monodromy
transformation of the singularity $C$, i.e., of the function germ
$f:(\C^2, 0)\to (\C, 0)$ such that $C=\{f=0\}$ (see, e.g., {\cite{A'C}}).
For $r >1$, one has $\zeta_C(t)=\Delta^C(t, \ldots, t)$ (for $r=1$,
$\zeta_C(t)=\Delta^C(t)/(1-t)$\,).

It was shown that all the coefficients of the Alexander polynomial
$\Delta^C(\tt)$ ($\tt=(t_1, \ldots, t_r)$\,)can be described as
Euler characteristics
of some spaces~-- complements to arrangements of projective hyperplanes
in projective spaces ({\cite{GDC}}). For a hypersurface singularity
of any dimension, the Lefschetz numbers of iterates of the classical
monodromy transformation have been described as Euler characteristics
of some subspaces in the space of (truncated) arcs ({\cite{DL}}).
The last result is connected with the theory of integration with
respect to the Euler characteristic in the space of arcs. Here we
discuss a similar notion (integration with respect to the Euler
characteristic) in the projectivization $\P{\cal O}_{\C^n,0}$
of the ring ${\cal O}_{\C^n,0}$  of germs of functions on $C^n$ at the origin
(here we consider it as a linear space)
and show that the Alexander polynomial and the zeta-function of a plane
curve singularity can be expressed as certain integrals over 
$\P{\cal O}_{\C^2,0}$ with respect to the Euler characteristic.

Let $J^k_{\C^n,0}$ be the space of $k$-jets of functions at the origin
in $(\C^n,0)$
($J^k_{\C^n,0}={\cal O}_{\C^n,0}/m^{k+1}\cong \C^{{n+k\choose k}}$,
where $m$ is the maximal ideal in ${\cal O}_{\C^n,0}$). For a complex linear
space $L$ (finite or infinite dimensional) let $\P L=(L\setminus\{0\})/\C^*$
be its projectivization, let $\P^* L$ be the disjoint union of $\P L$
with a point (in some sense $\P^* L=L/\C^*$). One has natural maps
$\pi_k: \P{\cal O}_{\C^n,0} \to \P^* J^k_{\C^n,0}$ and 
$\pi_{k,\ell}: \P^* J^k_{\C^n,0} \to \P^* J^\ell_{\C^n,0}$ for $k \ge \ell$.
Over $\P J^\ell_{\C^n,0} \subset \P^* J^\ell_{\C^n,0}$ the map $\pi_{k,\ell}$
is a locally trivial (and in fact trivial) fibration, the fibre of which
is a complex linear space of some dimension.

\begin{definition}
A subset $X\subset \P{\cal O}_{\C^n,0}$ is said to be cylindric if
$X=\pi_k^{-1}(Y)$ for a semi-algebraic subset
$Y\subset \P J^k_{\C^n,0} \subset \P^* J^k_{\C^n,0}$.
\end{definition}

\begin{definition}
For a cylinder subset $X\subset \P{\cal O}_{\C^n,0}$ ($X=\pi_k^{-1}(Y)$,
$Y\subset \P J^k_{\C^n,0}$) its Euler characteristic $\chi(X)$ is defined
as the Euler characteristic $\chi(Y)$ of the set $Y$.
\end{definition}

\begin{remark}
A semi-algebraic subset of a finite dimensional projective space
(e.g., the set $Y$ above) can be represented as the union of a finite number
of cells which do not intersect each other. The Euler characteristic
of such a set is defined as the alternative sum of numbers of cells
of different dimensions. Defined this way, the Euler characteristic
satisfies the additivity property:
$$
\chi(Y_1\cup Y_2)= \chi(Y_1) + \chi(Y_2) - \chi(Y_1\cap Y_2), 
$$
and therefore can be considered as a generalized (non-positive)
measure on the algebra of semi-algebraic subsets.
\end{remark}

Let $\psi: \P{\cal O}_{\C^n,0} \to A$ be a function with values in
an Abelian group $G$.

\begin{definition}
We say that the function $\psi$ is cylindric if, for each $a\ne 0$
the set $\psi^{-1}(a)\subset \P{\cal O}_{\C^n,0}$ is cylindric.
\end{definition}

\begin{definition}
The integral of a cylindric function $\psi$ over $\P{\cal O}_{\C^n,0}$
with respect to the Euler characteristic is
$$
\int_{\P{\cal O}_{\C^n,0}}\psi d\chi =
\sum_{a\in A, a\ne 0} \chi(\psi^{-1}(a))\cdot a
$$
if this sum has sense in $A$. If the integral exists (has sense)
the function $\psi$ is said to be integrable.
\end{definition}

\begin{remark}
In a similar way one can define a generalized Euler characteristic $[X]$
of a cylindric subset of $\P{\cal O}_{\C^n,0}$ (or of ${\cal O}_{\C^n,0}$)
with values in the Grothendieck ring of complex algebraic varieties
localized by the the class $\L$ of the complex line and thus the
corresponding notion of integration (see, e.g., {\cite{C}}). For that one
can define $[X]$ as $[Y]\cdot\L^{-{n+k\choose k}}$.
\end{remark}

For a plane curve singularity $C=\bigcup\limits_{i=1}^{r}C_i$, let
$\varphi_i:(\C_i, 0)\to(\C^2, 0)$ be parameterizations (uniformizations)
of the branches $C_i$ of the curve $C$ (i.e., ${\rm{Im}}\,\varphi_i=C_i$
and $\varphi_i$ is an isomorphism between $\C_i$ and $C_i$
outside of the origin). For a germ $g\in{\cal O}_{\C^2, 0}$,
let $v_i(g)$ be the power of the leading term in the power series
decomposition of the germ $g\circ\varphi_i:(\C_i,0)\to \C$:
$g\circ\varphi_i(\tau_i)=c_i\cdot \tau_i^{v_i}+{~terms~of~higher~degree~}$,
where $c_i\ne 0$. If $g\circ\varphi_i(t)\equiv 0$, $v_i(g)$ is assumed
to be equal to $\infty$. Let $\vv(g)=(v_1(g), \ldots, v_r(g))\in \Z_{\ge 0}^r$,
$v(g)= \Vert \vv(g) \Vert = v_1(g)+\ldots+v_r(g)$. Let $\Z[[t]]$ (respectively
$\Z[[t_1, \ldots, t_r]]$) be the group (with respect to addition) of formal
power series in the variable $t$ (respectively in $t_1$, \dots, $t_r$.
For $\vv=(v_1, \ldots, v_r)\in Z_{\ge 0}^r$, let
$\tt^\vv=t_1^{v^1}\cdot\ldots\cdot t_r^{v_r}$; we assume $t^\infty=0$.

\begin{theorem}
For each $\vv \in Z_{\ge 0}^r$,
the set $\{g\in \P{\cal O}_{\C^2, 0}: \vv(g)=\vv\}$ is cylindric. Therefore
the functions $\tt^{\vv(g)}$ and $t^{v(g)}$ on $\P{\cal O}_{\C^2, 0}$ with
values in $\Z[[t_1, \ldots, t_r]]$ and $\Z[[t]]$ respectively are cylindric.
\end{theorem}

{\bf Proof} follows from the fact that, for $g\in m^s$, $v_i(g)\ge s$,
i.e., the Taylor series of $g\circ\varphi_i(\tau_i)$ starts from terms
of degree at least $s$. Therefore the functions $\tt^{\vv(g)}$ and $t^{v(g)}$
on $\P{\cal O}_{\C^2, 0}$ are integrable (since
$\sum\limits_{\vv\in\Z_{\ge 0}^r}\ell(\vv)\tt^\vv\in \Z[[t_1, \cdots, t_r]]$
for any integers $\ell(\vv)$).

\begin{theorem}
For $r>1$, $$\int\limits_{\P{\cal O}_{\C^2, 0}}\tt^{\vv(g)} d\chi =
\Delta^C(t_1,\ldots, t_r);$$
for $r\ge 1$, $$\int\limits_{\P{\cal O}_{\C^2, 0}}t^{v(g)} d\chi =
\zeta_C(t).$$
\end{theorem}

{\bf Proof} follows from the results of {\cite{GDC}}.

\end{document}